%% file: arxiv.final.tex
\documentclass[10pt,twoside,reqno,draft]{amsart}
\input{prembu}
\usepackage{amsmath}
\usepackage{amssymb}
\usepackage{amsfonts}
\usepackage{amsthm}
\begin{document}
\setcounter{page}{1}

\leftline{{\em Accepted for Publication in}}
\leftline{\bf\em Fixed Point Theory}
\leftline{\em Cluj-Napoca, Romania} \leftline{\footnotesize
http://www.math.ubbcluj.ro/$^{\sim}$nodeacj/sfptcj.html}

\vs*{1.5cm}

\title[Nonlinear contractions in metric spaces]{\large Nonlinear contractions in metric spaces under
locally $T$-transitive binary relations}
\author[Aftab Alam and Mohammad Imdad]{Aftab Alam$^{*}$ and Mohammad Imdad$^{**}$}
\date{}
\maketitle

\vs*{-0.5cm}

\bc
{\footnotesize
$^{*}$Department of Mathematics, Aligarh Muslim University, Aligarh-202002, India\\
E-mail: aafu.amu@gmail.com\\
\medskip
$^{**}$Department of Mathematics, Aligarh Muslim University, Aligarh-202002, India\\
E-mail: mhimdad@gmail.com\\
}
\ec

\bigskip

{\footnotesize
\noindent
{\bf Abstract.} In this paper, we present
a variant of Boyd-Wong fixed point theorem in a metric space
equipped with a locally $T$-transitive binary relation,
which under universal relation reduces to Boyd-Wong (Proc. Amer. Math. Soc. {\bf 20}(1969), 458-464) and Jotic (Indian J. Pure Appl. Math. {\bf 26}(1995), 947-952) fixed point theorems. Also, our results extend
several other well-known fixed point theorems such as: Alam and
Imdad (J. Fixed Point Theory Appl. {\bf 17}(2015), no 4, 693-702) and Karapinar and Rold$\acute{\rm
a}$n-L$\acute{\rm o}$pez-de-Hierro (J. Inequal. Appl. {\bf 2014:522}(2014), 12 pages) besides some others.

\noindent
{\bf Key Words and Phrases}: $\mathcal{R}$-continuity; locally
$T$-transitive binary relations; $\varphi$-contractions; $\mathcal{R}$-connected sets.

\noindent {\bf 2010 Mathematics Subject Classification}: 47H10,
54H25.

}

\bigskip

\section{Introduction}
A variety of generalizations of the classical Banach contraction
principle \cite{B} is available in the existing literature of metric
fixed
point theory. These generalizations are obtained in the different directions such as:\\

\indent\hspace{0.5cm} (1) enlarging the class of ambient spaces,\\
\indent\hspace{0.5cm} (2) improving the underlying contraction condition,\\
\indent\hspace{0.5cm} (3) weakening the involved metrical notions
($e.g.$ completeness, continuity $etc.$).\\

Recently, Alam and Imdad \cite{RT1,RT2} obtained an interesting
generalization of classical Banach contraction principle by using an
amorphous (arbitrary) binary relation. In doing so, the authors
introduced the relation-theoretic analogues of certain involved metrical notions such as: contraction, completeness, continuity $etc.$.
In fact, under the universal relation, such
newly defined notions reduce to their corresponding usual notions and henceforth relation-theoretic metrical fixed/coincidence point
theorems reduce to their corresponding classical fixed/coincidence
point theorems (under the universal relation).\\

Recall that given a nonempty set $X$, a subset $\mathcal{R}$ of
$X^2$ is called a binary relation on $X$. For simplicity, we
sometimes write $x\mathcal{R}y$ instead of $(x,y)\in\mathcal{R}$.
Given $E\subseteq X$ and a binary relation $\mathcal{R}$ on $X$, the
restriction of $\mathcal{R}$ to $E$, denoted by $\mathcal{R}|_E$, is
defined to be the set $\mathcal{R}\cap E^2$
($i.e.~\mathcal{R}|_E:=\mathcal{R}\cap E^2$). Indeed,
$\mathcal{R}|_E$ is a relation on $E$ induced by $\mathcal{R}$.\\

Out of various classes of binary relations in practice, the
following ones are relevant in the present context.\\

A binary relation $\mathcal{R}$ on a nonempty set $X$ is called\\
\indent\hspace{0.5cm}$\bullet$ amorphous if $\mathcal{R}$ has no specific property at all,\\
\indent\hspace{0.5cm}$\bullet$ universal if $\mathcal{R}=X^2$,\\
\indent\hspace{0.5cm}$\bullet$ empty if $\mathcal{R}=\emptyset$,\\
\indent\hspace{0.5cm}$\bullet$ reflexive if $(x,x)\in \mathcal{R}~~~\forall x\in X$,\\
\indent\hspace{0.5cm}$\bullet$ symmetric if whenever $(x,y)\in \mathcal{R}$ then $(y,x)\in \mathcal{R}$,\\
\indent\hspace{0.5cm}$\bullet$ antisymmetric if whenever $(x,y)\in \mathcal{R}$ and $(y,x)\in \mathcal{R}$ then $x=y$,\\
\indent\hspace{0.5cm}$\bullet$ transitive if whenever $(x,y)\in \mathcal{R}$ and $(y,z)\in \mathcal{R}$ then $(x,z)\in \mathcal{R}$,\\
\indent\hspace{0.5cm}$\bullet$ complete if $(x,y)\in \mathcal{R}~or~(y,x)\in \mathcal{R}~~\forall~ x,y\in X$,\\
\indent\hspace{0.5cm}$\bullet$ partial order if $\mathcal{R}$ is reflexive, antisymmetric and transitive.\\

Throughout this paper, $\mathcal{R}$ stands for a nonempty binary
relation but for the sake of simplicity, we write only `binary
relation' instead of `nonempty binary relation'. Also, $\mathbb{N}$
stands for the set of natural numbers, while $\mathbb{N}_0$ for the
set of whole numbers ($i.e.~\mathbb{N}_0:=\mathbb{N} \cup
\{0\}$).\\

The following family of control functions is introduced by Lakshmikantham and \'{C}iri\'{c} \cite{C2}.\\
$$\Phi=\Big\{\varphi:[0,\infty)\to [0,\infty):\varphi(t)<t\;{\rm for~each}\; t>0\;{\rm and~}\lim\limits_{r\to t^+}\varphi(r)<t \;{\rm
for~each~t>0}\Big\}.$$

With a view to have a self-contained presentation, we recall two fixed point theorems involving nonlinear contractions for class $\Phi$ using partial order and transitive binary relations respectively, which have inspired our results in the present paper. In order to understand such results better, we recall firstly the relevant definitions and thereafter state the corresponding results.\\

\noindent{\bf Definition 1.1} \cite{P14,P7,P8}. Let $X$ be a
nonempty set equipped with a partial order $\preceq$. A self-mapping
$T$ on $X$ is called increasing or isotone or
order-preserving if for any $x,y\in X$,
$$x\preceq y \Rightarrow T(x)\preceq T(y).$$

The following notion is formulated by using a suitable property with a view to avoid the necessity of the continuity requirement on the involved mapping specially in the
hypotheses of a fixed point theorem due to Nieto and
Rodr\'{\i}guez-L\'{o}pez \cite{PF2}.\\

\noindent{\bf Definition 1.2} \cite{PGF13}. Let $(X,d)$ be a
metric space equipped with a partial order $\preceq$. We say that
$(X,d,\preceq)$ has {\it ICU}\;(increasing-convergence-upper bound)
property if every increasing convergent sequence in $X$ is bounded
above by its limit (as an upper bound).\\

The following result, a variant of fixed point theorem of Nieto and
Rodr\'{\i}guez-L\'{o}pez \cite{PF2} under $\varphi$-contraction, is
contained in many papers ( $e.g.$ Wu and Liu
(\cite{PGF12}, Theorem 2.1), Samet $et\;al.$ (\cite{R1}, Remark
1.3), Kutbi $et\;al.$ (\cite{FIC8}, Theorem 5), Karapinar $et\;al.$
(\cite{FIC9}, Theorem 10) and Karapinar and Rold$\acute{\rm
a}$n-L$\acute{\rm o}$pez-de-Hierro (\cite{FIC10},
Theorem 1.2)).\\

\noindent{\bf Theorem 1.3}. {\it Let $(X,d)$ be a metric space equipped with a partial order $\preceq$ and $T$ a self-mapping on $X$. Suppose that the following conditions hold:\\
\indent\hspace{0.5mm}$(a)$ $(X,d)$ is complete,\\
\indent\hspace{0.5mm}$(b)$ $T$ is increasing,\\
\indent\hspace{0.5mm}$(c)$ either $T$ is continuous or $(X,d,\preceq)$ has {\it ICU} property,\\
\indent\hspace{0.5mm}$(d)$ there exists $x_{0}\in X$ such that
$x_{0} \preceq
T(x_{0})$,\\
\indent\hspace{0.5mm}$(e)$ there exists $\varphi\in \Phi$ such that \\
\indent\hspace{2.5cm}$d(Tx,Ty)\leq\varphi(d(x,y))\;\;\forall ~x,y\in X$ with $x\preceq y$.\\
Then $T$ has a fixed point. Moreover, if for all $x,y\in X$, there
exists $z\in X$ such that $x\preceq z$ and $y\preceq z$, then we obtain
uniqueness of the fixed point.}\\

The above result seems natural but the partial order relation is
very restrictive. Samet and Vetro \cite{NFI} introduced the notion
of $F$-invariant set and utilized the same to
prove some coupled fixed point results for generalized linear
contractions in metric spaces without partial order. In 2012,
Sintunaravat $et\;al.$ \cite{FIC3} introduced the notion of
transitive property and utilized the same to extend some
Samet-Vetro coupled fixed point theorems for nonlinear
contractions. On the other hand, Kutbi $et\;al.$ \cite{FIC8} weakened the
notion of $F$-invariant sets by introducing the notion of $F$-closed
sets. Recently, Karapinar $et\;al.$ \cite{FIC9} proved some
unidimensional versions of earlier coupled fixed point results
involving $F$-closed sets and then obtained such coupled fixed point
results as easy by using their corresponding (unidimensional) fixed point results. As noticed in Alam and
Imdad \cite{RT2}, the relation-theoretic metrical
fixed/coincidence point theorems combine the idea contained in
Karapinar $et\;al.$ \cite{FIC9} as the set {\it M} (utilized by
Karapinar $et\;al.$ \cite{FIC9}) being subset of $X^2$ is in fact a binary relation on $X$.\\

The following notions are unidimensional formulations of transitive
property and $F$-closed sets.\\

\noindent{\bf Definition 1.4} \cite{FIC9, FIC11}. We say that a
nonempty subset {\it M}$\subseteq X^2$ is\\
\indent\hspace{0.5cm}$\bullet$ transitive if $(x,z)\in {\it M}$ for all $x,y,z\in X$ such that $(x,y),(y,z)\in {\it M}$.\\
Given a mapping $T:X\rightarrow X$, we say that {\it M} is\\
\indent\hspace{0.5cm}$\bullet$ $T$-transitive if $(Tx,Ty)\in{\it M}$ for all $x,y,z\in X$ such that $(Tx,Ty),(Ty,Tz)\in {\it M}$,\\
\indent\hspace{0.5cm}$\bullet$ $T$-closed if $(Tx,Ty)\in{\it M}$ for all $x,y\in X$ such that $(x,y)\in {\it M}$.\\

\noindent{\bf Definition 1.5} \cite{FIC9}. Let $(X, d)$ be a metric
space and let {\it M}$\subseteq X^2$ be a subset. We say that
$(X,d,{\it M})$ is regular if for all sequence $\{x_n\}\subseteq X$
such that $x_n\stackrel{d}{\longrightarrow} x$ and $(x_n, x_{n+1})\in {\it M}$ for all
$n\in \mathbb{N}$, we have $(x_n,x)\in {\it M}$ for all $n\in \mathbb{N}$.\\

The following fixed point theorem indicated in Karapinar and
Rold$\acute{\rm a}$n-L$\acute{\rm o}$pez-de-Hierro \cite{FIC10} is
a unidimensional version of coupled fixed point
theorem of Sintunaravat $et\;al.$ \cite{FIC3}.\\

\noindent{\bf Theorem 1.6} \cite{FIC10}. {\it Let $(X,d)$ be a metric
space, let $T:X\rightarrow X$ be a mapping and let {\it M}$\subseteq
X^2$ be a subset such that\\
\indent\hspace{0.5mm}$(a)$ $(X,d)$ is complete,\\
\indent\hspace{0.5mm}$(b)$ {\it M} is $T$-closed and transitive,\\
\indent\hspace{0.5mm}$(c)$ either $T$ is continuous or $(X,d,{\it
M})$ is regular,\\
\indent\hspace{0.5mm}$(d)$ there exists $x_{0}\in X$ such that $(x_{0},Tx_{0})\in {\it M}$,\\
\indent\hspace{0.5mm}$(e)$ there exists $\varphi\in \Phi$ such that \\
\indent\hspace{2.5cm}$d(Tx,Ty)\leq\varphi(d(x,y))\;\;\forall~ x,y\in X$ with $(x,y)\in {\it M}$.\\

Then $T$ has, at least, a fixed point.}\\

The aim of this paper is to extend the Alam-Imdad relation-theoretic fixed
point theorem \cite{RT1} for nonlinear contractions. Our results improve Theorems 1.3 and 1.6 in the following respects:
\begin{enumerate}
\item [{(i)}] the underlying binary relations (partial order or transitive) are replaced by an optimal condition of transitivity namely: locally $T$-transitive binary relation,
\item [{(ii)}] the nonlinear contractive class $\Phi$ is replaced by relatively enlarger class due to Boyd and Wong \cite{B2},
\item [{(iii)}] the involved metrical notions namely: completeness and continuity are replaced by their $\mathcal{R}$-analogues,
\item [{(iv)}] the {\it ICU} property and the regularity of $X$ are replaced by relatively weaker notion namely: $d$-self-closedness.
\end{enumerate}

\section{Relation-theoretic notions and auxiliary results}
In this section, for the sake of completeness, we summarize some
necessary definitions and basic results related to our main results.\\

\noindent{\bf Definition 2.1} \cite{RT1}. Let $\mathcal{R}$ be a
binary relation on a nonempty set $X$ and $x,y\in X$. We say that
$x$ and $y$ are $\mathcal{R}$-comparative if either $(x,y)\in
\mathcal{R}$ or $(y,x)\in \mathcal{R}$. We denote it by $[x,y]\in
\mathcal{R}$.\\

\noindent{\bf Definition 2.2} \cite{ST}. Let $X$ be a nonempty set
and $\mathcal{R}$ a binary relation on $X$.
\begin{enumerate}
\item [{(1)}] The inverse or transpose or dual relation of $\mathcal{R}$, denoted by $\mathcal{R}^{-1}$,
 is
defined by $\mathcal{R}^{-1}=\{(x,y)\in X^2:(y,x)\in \mathcal{R}\}$.
\item [{(2)}] The symmetric closure of $\mathcal{R}$, denoted by $\mathcal{R}^s$,
is defined to be the set $\mathcal{R}\cup \mathcal{R}^{-1}$
($i.e.~\mathcal{R}^s:=\mathcal{R}\cup \mathcal{R}^{-1}$). Indeed,
$\mathcal{R}^s$ is the smallest symmetric relation on $X$ containing
$\mathcal{R}$.
\end{enumerate}

\noindent{\bf Proposition 2.3} \cite{RT1}. {\it For a binary relation
$\mathcal{R}$ defined on a nonempty set $X$,
$$(x,y)\in\mathcal{R}^s\Longleftrightarrow [x,y]\in\mathcal{R}.$$}

\noindent{\bf Definition 2.4} \cite{RT1}. Let $X$ be a nonempty
set and $\mathcal{R}$ a binary relation on $X$. A sequence
$\{x_n\}\subset X$ is called $\mathcal{R}$-preserving if
$$(x_n,x_{n+1})\in\mathcal{R}\;\;\forall~n\in \mathbb{N}_{0}.$$

\noindent{\bf Definition 2.5} \cite{RT1}. Let $X$ be a nonempty
set and $T$ a self-mapping on $X$. A binary relation $\mathcal{R}$ on $X$ is called $T$-closed if for any $x,y\in X$,
$$(x,y)\in \mathcal{R}\Rightarrow (Tx,Ty)\in \mathcal{R}.$$

\noindent{\bf Proposition 2.6} \cite{RT1}. {\it Let $X$ be a nonempty
set, $\mathcal{R}$ a binary relation on $X$ and $T$ a self-mapping
on $X$. If $\mathcal{R}$ is $T$-closed, then
$\mathcal{R}^s$ is also $T$-closed.}\\

\noindent{\bf Proposition 2.7}. {\it Let $X$ be a nonempty set,
$\mathcal{R}$ a binary relation on $X$ and $T$ a self-mapping on
$X$. If $\mathcal{R}$ is $T$-closed, then, for all $n\in
\mathbb{N}_0$, $\mathcal{R}$ is also $T^n$-closed, where $T^n$ denotes $n$th iterate of $T$.}\\

\noindent{\bf Definition 2.8} \cite{RT2}. Let $(X,d)$ be a metric
space and $\mathcal{R}$ a binary relation on $X$. We say that
$(X,d)$ is $\mathcal{R}$-complete if every $\mathcal{R}$-preserving
Cauchy sequence in $X$ converges.\\

Clearly, every complete metric space is $\mathcal{R}$-complete, for
any binary relation $\mathcal{R}$. Particularly, under the universal
relation the notion of
$\mathcal{R}$-completeness coincides with usual completeness.\\

\noindent{\bf Definition 2.9} \cite{RT2}. Let $(X,d)$ be a metric
space, $\mathcal{R}$ a binary relation on $X$ and $x\in X$. A
self-mapping $T$ on $X$ is called $\mathcal{R}$-continuous at $x$
if for any $\mathcal{R}$-preserving sequence $\{x_n\}$ such that
$x_n\stackrel{d}{\longrightarrow} x$, we have
$T(x_n)\stackrel{d}{\longrightarrow} T(x)$. Moreover, $T$ is called
$\mathcal{R}$-continuous if it is $\mathcal{R}$-continuous at each
point of $X$.\\

Clearly, every continuous mapping is $\mathcal{R}$-continuous, for
any binary relation $\mathcal{R}$. Particularly, under the universal
relation the notion of $\mathcal{R}$-continuity coincides with usual
continuity.\\

The following notion is a generalization of $d$-self-closedness of a
partial order relation $(\preceq)$ (defined by Turinici
\cite{T-PFP,T-LCP}).\\

\noindent{\bf Definition 2.10} \cite{RT1}. Let $(X,d)$ be a metric
space. A binary relation $\mathcal{R}$ on $X$ is called
$d$-self-closed if for any $\mathcal{R}$-preserving sequence
$\{x_n\}$ such that $x_n\stackrel{d}{\longrightarrow} x$, there
exists a subsequence $\{x_{n_k}\}{\rm \;of\;} \{x_n\} \;{\rm
with}\;\;[x_{n_k},x]\in\mathcal{R}~~~\forall~k\in \mathbb{N}_{0}.$\\

\noindent{\bf Definition 2.11} \cite{BR1}. Let $X$ be a nonempty
set and $\mathcal{R}$ a binary relation on $X$. A subset $E$ of $X$
is called $\mathcal{R}$-directed if for each pair $x,y\in E$, there
exists $z\in X$ such that $(x,z)\in\mathcal{R}$ and
$(y,z)\in\mathcal{R}$.\\

\noindent{\bf Definition 2.12} \cite{DM}. Let $X$ be a nonempty set
and $\mathcal{R}$ a binary relation on $X$. For $x,y\in X$, a path
of length $k$ (where $k$ is a natural number) in $\mathcal{R}$ from
$x$ to $y$ is a finite sequence
$\{z_0,z_1,z_2,...,z_{k}\}\subset X$ satisfying the following conditions:\\
\indent\hspace{0.5mm} (i) $z_0=x~{\rm and}~z_k=y$,\\
\indent\hspace{0.5mm} (ii) $(z_i,z_{i+1})\in\mathcal{R}$ for each $i~(0\leq i\leq k-1)$.\\
Notice that a path of length $k$ involves $k+1$ elements of $X$,
although they are not necessarily distinct.\\

\noindent{\bf Definition 2.13} \cite{RT2}. Let $X$ be a nonempty
set and $\mathcal{R}$ a binary relation on $X$. A subset $E$ of $X$
is called $\mathcal{R}$-connected if for each pair $x,y\in E$, there
exists a path in $\mathcal{R}$ from $x$ to $y$.\\

\noindent Given a binary relation $\mathcal{R}$ and a self-mapping
$T$ on a nonempty set $X$, we use the following notations.\\

\indent\hspace{0.5mm} (i) $F(T)$:=the set of all fixed points of $T$,\\
\indent\hspace{0.5mm} (ii) $X(T,\mathcal{R}):=\{x\in X:(x,Tx)\in \mathcal{R}\}$.\\

The following result is the relation-theoretic version of Banach contraction principle.\\

\noindent{\bf Theorem 2.14} \cite{RT1,RT2}. {\it Let $(X,d)$ be a metric
space, $\mathcal{R}$ a binary relation on $X$ and $T$ a self-mapping
on $X$. Suppose that the following
conditions hold:\\
\indent\hspace{0.5mm}$(a)$ $(X,d)$ is $\mathcal{R}$-complete,\\
\indent\hspace{0.5mm}$(b)$ $\mathcal{R}$ is $T$-closed,\\
\indent\hspace{0.5mm}$(c)$ either $T$ is $\mathcal{R}$-continuous or $\mathcal{R}$ is $d$-self-closed,\\
\indent\hspace{0.5mm}$(d)$ $X(T,\mathcal{R})$ is nonempty,\\
\indent\hspace{0.5mm}$(e)$ there exists $\alpha\in [0,1)$ such that \\
\indent\hspace{2.5cm}$d(Tx,Ty)\leq\alpha d(x,y)\;\;\forall~ x,y\in X$ with $(x,y)\in \mathcal{R}$.\\
Then $T$ has a fixed point. Moreover if $X$ is $\mathcal{R}^s$-connected then $T$ has a unique fixed point.}\\

Now, we re-define the notion of $T$-transitivity employed in
Definition 1.4 in the framework of binary relation.\\

\noindent{\bf Definition 2.15}. Let $X$ be a nonempty
set and $T$ a self-mapping on $X$. A binary relation $\mathcal{R}$ on $X$ is called
$T$-transitive if for any $x,y,z\in X$,
$$(Tx,Ty),(Ty,Tz)\in\mathcal{R}\Rightarrow(Tx,Tz)\in\mathcal{R}.$$

Inspired by Turinici \cite{BR6,BR7}, we introduce the following notion by localizing the notion of transitivity.\\

\noindent{\bf Definition 2.16}. A binary relation
$\mathcal{R}$ on a nonempty set $X$ is called locally
transitive if for each (effectively) $\mathcal{R}$-preserving sequence $\{x_n\}\subset X$ (with range
$E:=\{x_n:n \in \mathbb{N}_{0}\}$), the binary relation $\mathcal{R}|_E$ is transitive.\\

Henceforth, the notions of $T$-transitivity and locally
transitivity both are relatively weaker than the notion of
transitivity but they are independent of each others. In order to
make them compatible, we introduce the following notion of
transitivity.\\

\noindent{\bf Definition 2.17}. Let $X$ be a nonempty
set and $T$ a self-mapping on $X$. A binary relation $\mathcal{R}$ on $X$ is called locally $T$-transitive if for each (effectively) $\mathcal{R}$-preserving sequence $\{x_n\}\subset T(X)$ (with range
$E:=\{x_n:n \in \mathbb{N}_{0}\}$), the binary relation $\mathcal{R}|_E$ is transitive.\\

The following result establishes the superiority of locally
$T$-transitivity over other types of transitivity.\\

\noindent{\bf Proposition 2.18}. {\it Let $X$ be a nonempty set,
$\mathcal{R}$ a binary relation on $X$ and $T$ a self-mapping on
$X$.
\begin{enumerate}
\item[{(i)}] $\mathcal{R}$ is $T$-transitive
$\Leftrightarrow$ $\mathcal{R}|_{T(X)}$ is transitive.
\item[{(ii)}] $\mathcal{R}$ is locally $T$-transitive
$\Leftrightarrow$ $\mathcal{R}|_{T(X)}$ is locally
transitive.
\item[{(iii)}] $\mathcal{R}$ is transitive $\Rightarrow$ $\mathcal{R}$ is locally
transitive $\Rightarrow$ $\mathcal{R}$ is locally
$T$-transitive.
\item[{(iv)}] $\mathcal{R}$ is transitive $\Rightarrow$ $\mathcal{R}$ is $T$-transitive $\Rightarrow$ $\mathcal{R}$ is locally $T$-transitive.\\
\end{enumerate}}
The following family of control functions is indicated in  Boyd and
Wong \cite{B2} but was later used in Jotic \cite{B4}.\\
$$\Omega=\Big\{\varphi:[0,\infty)\to [0,\infty):\varphi(t)<t\;{\rm for~each}\; t>0\;{\rm and~}\limsup\limits_{r\to t^+}\varphi(r)<t \;{\rm
for~each~t>0}\Big\}.$$

It is clear that the class $\Omega$ enlarges the class $\Phi$, $i.e.$, $\Phi \subset \Omega$.\\

\noindent{\bf Proposition 2.19}. {\it If $(X,d)$ is a metric space,
$\mathcal{R}$ is a binary relation on $X$, $T$ is a
self-mapping on $X$ and $\varphi\in \Omega$, then the following contractivity conditions are equivalent:\\
\indent\hspace{0.1cm}(I) $d(Tx,Ty)\leq\varphi(d(x,y))\;\;\forall~ x,y\in X$ with $(x,y)\in \mathcal{R}$,\\
\indent\hspace{0.1cm}(II) $d(Tx,Ty)\leq\varphi(d(x,y))\;\;\forall~ x,y\in X$ with $[x,y]\in \mathcal{R}$.}\\

We skip the proof of above proposition as it is similar to that of
Proposition 2.3 \cite{RT1}.\\

Finally, we record the following known results, which are needed in the proof of our main results. \\

\noindent{\bf Lemma 2.20} \cite{PGF13}. {\it Let $\varphi\in\Omega$. If
$\{a_n\}\subset (0,\infty)$ is a sequence such that $a_{n+1}\leq
\varphi(a_n)\; \forall~ n\in \mathbb{N}_0$, then
$\lim\limits_{n\to\infty}a_{n}=0$.}\\

\noindent{\bf Lemma 2.21} \cite{BR5,BR7}. {\it Let $(X,d)$ be a metric space and
$\{x_n\}$ a sequence in $X$. If
$\{x_n\}$ is not a Cauchy, then there exist $\epsilon>0$ and two subsequences $\{x_{n_k}\}$ and $\{x_{m_k}\}$ of $\{x_n\}$ such that\\
\indent\hspace{0.5mm}(i) $k\leq m_k<n_k\;\;\forall~k\in\mathbb{N}$,\\
\indent\hspace{0.5mm}(ii) $d(x_{m_k},x_{n_k})> \epsilon\;\;\forall~k\in\mathbb{N}$,\\
\indent\hspace{0.5mm}(iii) $d(x_{m_k},x_{n_{k-1}})\leq\epsilon\;\;\forall~k\in\mathbb{N}$.\\
Moreover, suppose that $\lim\limits_{n\to
\infty}d(x_n,x_{n+1})=0$, then\\
\indent\hspace{0.5mm}(iv) $\lim\limits_{k\to
\infty}d(x_{m_k},x_{n_k})=\epsilon,$\\
\indent\hspace{0.5mm}(v) $\lim\limits_{k\to \infty}
d(x_{m_k+1},x_{n_k+1})=\epsilon.$}

\section{Fixed Point Theorems}

Firstly, we prove a result on the existence of fixed points under $\varphi$-contractivity condition, which runs as follows.\\

\noindent{\bf Theorem 3.1}. {\it Let $(X,d)$ be a metric space,
$\mathcal{R}$ a binary relation on $X$ and $T$ a self-mapping on
$X$. Suppose that the following
conditions hold:\\
\indent\hspace{0.5mm}$(a)$ $(X,d)$ is $\mathcal{R}$-complete,\\
\indent\hspace{0.5mm}$(b)$ $\mathcal{R}$ is $T$-closed and locally $T$-transitive,\\
\indent\hspace{0.5mm}$(c)$ either $T$ is $\mathcal{R}$-continuous or $\mathcal{R}$ is $d$-self-closed,\\
\indent\hspace{0.5mm}$(d)$ $X(T,\mathcal{R})$ is nonempty,\\
\indent\hspace{0.5mm}$(e)$ there exists $\varphi\in \Omega$ such that \\
\indent\hspace{2.5cm}$d(Tx,Ty)\leq\varphi(d(x,y))\;\;\forall~ x,y\in X$ with $(x,y)\in \mathcal{R}$.\\
Then $T$ has a fixed point.}\\

\noindent{\bf Proof}. In view of assumption $(d)$, take arbitrarily
$x_0\in X(T,\mathcal{R})$. Construct the sequence $\{x_n\}$ of Picard
iterates based at the initial point $x_0$, $i.e$,
$$x_n=T^n(x_0)\;\forall~ n\in \mathbb{N}_0.\eqno(1)$$

As $(x_0,Tx_0)\in \mathcal{R}$, using $T$-closedness of
$\mathcal{R}$ and Proposition 2.7, we obtain
$$(T^nx_0,T^{n+1}x_0)\in
\mathcal{R}$$ so that $$(x_n,x_{n+1})\in \mathcal{R}\;\;\forall~n\in
\mathbb{N}_0.\eqno(2)$$ Thus the sequence $\{x_n\}$ is
$\mathcal{R}$-preserving. Applying the contractivity condition $(e)$
to (2), we deduce, for all $n\in \mathbb{N}
 _0$ that
$$d(x_{n+1},x_{n+2})\leq\varphi(d(x_{n},x_{n+1})).$$
Hence by Lemma 2.20, we obtain
$$\lim\limits_{n\to \infty}d(x_n, x_{n+1})=0.\eqno(3)$$
Now, we show that $\{x_n\}$ is a Cauchy sequence. On contrary,
suppose that $\{x_n\}$ is not Cauchy. Therefore, owing to Lemma 2.21,
there exist $\epsilon>0$ and two subsequences $\{x_{n_k}\}$ and
$\{x_{m_k}\}$ of $\{x_n\}$ such that

$$k\leq m_k<n_k,\;d(x_{m_k},x_{n_k})>\epsilon \geq d(x_{m_k},x_{n_{k-1}}) \;\;\forall~k\in\mathbb{N}.\eqno(4)$$ Further, in view of
(3), Lemma 2.21 assures us that
$$\lim_{k\rightarrow\infty}\
d(x_{m_k},x_{n_k})=\displaystyle\lim_{k\rightarrow\infty}\
d(x_{m_k+1},x_{n_k+1})=\epsilon.\eqno(5)$$

Denote $r_k:=d(x_{m_k},x_{n_k})$. As $\{x_n\}$ is
$\mathcal{R}$-preserving (owing to (2)) and $\{x_n\}\subset T(X)$ (owing to (1)), by locally
$T$-transitivity of $\mathcal{R}$, we have $(x_{m_k},x_{n_k})\in \mathcal{R}$. Hence, applying contractivity condition
$(e)$, we obtain
\begin{eqnarray*}
d(x_{m_k+1},x_{n_k+1})
&=&d(Tx_{m_k},Tx_{n_k})\\
&\leq&\varphi(d(x_{m_k},x_{n_k})).\\
&=&\varphi(r_k)
\end{eqnarray*}
so that $$d(x_{m_k+1},x_{n_k+1})\leq \varphi(r_k).\eqno(6)$$
Using the facts that $r_k\longrightarrow\epsilon$ in the real line as $k\rightarrow \infty$ (owing to (5)) and $r_k>\epsilon\;\forall~k\in\mathbb{N}$ (owing to (4)) and by the definition of $\Omega$, we have
$$\displaystyle\limsup_{k\rightarrow\infty}
\varphi(r_k)=\displaystyle\limsup_{r\rightarrow\epsilon^+}\
\varphi(r)<\epsilon.\eqno(7)$$
On taking limit superior as $k\longrightarrow\infty$ in (6) and using (5) and
(7), we obtain
$$\epsilon=\displaystyle\limsup_{k\rightarrow\infty}d(x_{m_k+1},x_{n_k+1})\leq\displaystyle\limsup_{k\rightarrow\infty}
\varphi(r_k)<\epsilon,$$ which is a contradiction. Therefore,
$\{x_n\}$ is a Cauchy sequence. Hence, $\{x_n\}$ is an
$\mathcal{R}$-preserving Cauchy sequence. By
$\mathcal{R}$-completeness of $X$, $\exists~x\in X$ such that
$x_n\stackrel{d}{\longrightarrow} x$.\\

Finally, we use assumption $(c)$ to show that $x$ is a fixed
point of $T$. Suppose that $T$ is $\mathcal{R}$-continuous. As
$\{x_n\}$ is $\mathcal{R}$-preserving with
$x_n\stackrel{d}{\longrightarrow} x$, $\mathcal{R}$-continuity of
$T$ implies that $x_{n+1}=T(x_n)\stackrel{d}{\longrightarrow} T(x)$.
Using the uniqueness of limit, we obtain $T(x)=x$, $i.e$, $x$ is
a fixed point of $T$.\\

Alternately, assume that $\mathcal{R}$ is $d$-self-closed. Again as
$\{x_n\}$ is $\mathcal{R}$-preserving such that
$x_n\stackrel{d}{\longrightarrow} x$, $d$-self-closedness of
$\mathcal{R}$ guarantees the existence of a subsequence
$\{x_{n_k}\}$ of $\{x_n\}$ with
$[x_{n_k},x]\in\mathcal{R}~~~\forall~k\in \mathbb{N} _0.$ On using the fact $[x_{n_k},x]\in\mathcal{R}$, assumption
$(e)$ and Proposition 2.19, we obtain
$$d(x_{n_k+1},Tx)=d(Tx_{n_k},Tx)\leq \varphi(d(x_{n_k},x))\;\;\forall~ k\in \mathbb{N}_0.$$
We claim that
$$d(x_{n_k+1},Tx)\leq d(x_{n_k},x)\;\;\forall~ k\in \mathbb{N}.\eqno(8)$$
On account of two different possibilities arising here, we consider
a partition $\{\mathbb{N}^0,\mathbb{N}^+\}$ of
$\mathbb{N},\;i.e.,\;\mathbb{N}^0\cup\mathbb{N}^+=\mathbb{N}\;{\rm
and~}\mathbb{N}^0\cap\mathbb{N}^+=\emptyset$ verifying that\\
\indent\hspace{1cm}(i) $d(x_{n_k},x)=0\;\;\forall~ k\in \mathbb{N}^0,$\\
\indent\hspace{1cm}(ii) $d(x_{n_k},x)>0\;\;\forall~ k\in \mathbb{N}^+.$\\
In case (i), we have $d(Tx_{n_k},Tx)=0\;\forall~ k\in \mathbb{N}^0$,
which implies that $d(x_{n_k+1},Tx)=0\;\forall~ k\in \mathbb{N}^0$
and hence (8) holds for all $k\in \mathbb{N}^0.$ In case (ii), by the definition of $\Omega$, we have $d(x_{n_k+1},Tx)\leq
\varphi(d(x_{n_k},x))<d(x_{n_k},x)\;\forall~ k\in \mathbb{N}^+$ and
hence (8) holds for all $k\in \mathbb{N}^+$. Thus
(8) holds for all $k\in \mathbb{N}.$\\
Taking limit of (8) as $k\rightarrow \infty$ and using
$x_{n_k}\stackrel{d}{\longrightarrow} x$, we obtain
$x_{n_k+1}\stackrel{d}{\longrightarrow} T(x)$. Owing to the
uniqueness of limit, we obtain $T(x)=x$ so that $x$ is a
fixed point of $T$. \\

Using Proposition 2.18, we obtain the following consequence of Theorem
3.1.\\

\noindent{\bf Corollary 3.2}. {\it Theorem 3.1 remains true if locally $T$-transitivity of $\mathcal{R}$ (utilized in assumption
$(b)$) is replaced by any one of the following
conditions (besides retaining rest of the hypotheses):\\
\indent\hspace{0.5mm}(i) $\mathcal{R}$ is transitive,\\
\indent\hspace{0.5mm}(ii) $\mathcal{R}$ is $T$-transitive,\\
\indent\hspace{0.5mm}(iii) $\mathcal{R}$ is locally transitive.}\\

Now, we prove a corresponding uniqueness result.\\

\noindent{\bf Theorem 3.3}. {\it In addition to the hypotheses of Theorem 3.1, suppose that the following condition holds:\\
\indent\hspace{0.5mm} $(u)$ $T(X)$ is
$\mathcal{R}^s$-connected.\\
Then $T$ has a unique fixed point.}\\

\noindent{\bf Proof}.  In view of Theorem 3.1, $F(T)\neq \emptyset$. Take
$x,y\in F(T)$, then for all $n \in \mathbb{N}_0$, we have
$$T^n(x)=x\;{\rm and}\;T^n(y)=y.\eqno(9)$$
Clearly $x,y\in T(X)$. By assumption $(u)$, there exists a path  (say
$\{z_0,z_1,z_2,...,z_k\}$) of some finite length $k$ in
$\mathcal{R}^s$ from $x$ to $y$ so that
$$z_0=x,\;z_{k}=y\;{\rm and}\;[z_i,z_{i+1}]\in\mathcal{R}\;{\rm for~each}\;i~(0\leq i\leq k-1).\eqno(10)$$
As $\mathcal{R}$ is $T$-closed, using Propositions 2.6 and 2.7, we
have
$$[T^nz_i,T^nz_{i+1}]\in\mathcal{R}\;{\rm for~each}\;i~(0\leq i\leq k-1)\;{\rm and}\;{\rm for~each}\;n\in
\mathbb{N}_0.\eqno(11)$$ Now, for each $n \in \mathbb{N}_0$ and for
each $i\;(0\leq i\leq k-1)$, define $t_n^i:=d(T^nz_i,T^nz_{i+1})$.
We claim that
$$\lim\limits_{n\to\infty}t_n^i=0.\eqno(12)$$
Fix $i$ and distinguish two cases. Firstly, suppose that
$t_{n_0}^i=d(T^{n_0}z_i,T^{n_0}z_{i+1})=0$ for some $n_0\in
\mathbb{N}_0$, $i.e.$ $T^{n_0}(z_i)=T^{n_0}(z_{i+1})$, which implies
that $T^{n_0+1}(z_i)=T^{n_0+1}(z_{i+1})$. Consequently, we get
$t_{n_0+1}^i=d(T^{n_0+1}z_i,T^{n_0+1}z_{i+1})=0$. Thus by induction,
we get $t_n^i=0\;\forall~ n\geq n_0,$ yielding thereby
$\lim\limits_{n\to\infty}t_n^i=0$. On the other hand, suppose that
$t_{n}>0\;\forall~n\in \mathbb{N}_0$, then on using (11), assumption
$(e)$ and Proposition 2.19, we obtain
\begin{eqnarray*}
t_{n+1}^i&=& d(T^{n+1}z_i,T^{n+1}z_{i+1})\\
&\leq& \varphi (d(T^nz_i,T^nz_{i+1}))\\
&=&\varphi(t_{n}^i)
\end{eqnarray*}
so that $$t_{n+1}^i\leq \varphi(t_{n}^i).$$ Hence, on applying Lemma
2.20, we obtain $\lim\limits_{n\to\infty}t_{n}^i=0$. Thus, in both the
cases, (12) is proved for each $i\;(0\leq i\leq k-1)$.

Making use of (9), (10), (12) and the triangular inequality, we obtain
$$d(x,y)= d(T^nz_0,T^nz_{k})\leq t_n^0+t_n^1+\cdots+t_n^{k-1}
\to 0\;\; as \;\; n\to\infty$$
so that $x=y$. Hence $T$ has a unique fixed point.\\

The following consequence of Theorem 3.3 is worth recording.\\

\noindent{\bf Corollary 3.4}. {\it Theorem 3.3 remains true if we replace
the condition $(u)$ by one of the following conditions (besides retaining rest of the hypotheses):\\
\indent\hspace{0.5mm} $(u^\prime)$ $\mathcal{R}|_{T(X)}$ is complete,\\
\indent\hspace{0.5mm} $(u^{\prime\prime})$ $T(X)$ is
$\mathcal{R}^s$-directed.}\\
\noindent{\bf Proof}. If $(u^\prime)$ holds, then for each $u,v\in T(X)$,
$[u,v]\in\mathcal{R}$, which
amounts to say that $\{u,v\}$ is a path of length 1 in
$\mathcal{R}^s$ from $u$ to $v$. Hence $T(X)$ is
$\mathcal{R}^s$-connected consequently
Theorem 3.3 gives rise the conclusion.\\
Otherwise, if $(u^{\prime\prime})$ holds then for each $u,v\in T(X)$,
$\exists~z\in X$ such that $[u,z]\in\mathcal{R}$ and
$[v,z]\in\mathcal{R}$, which
amounts to say that $\{u,z,v\}$ is a path of length 2
in $\mathcal{R}^s$ from $u$ to $v$. Hence
$T(X)$ is
$\mathcal{R}^s$-connected and again by
Theorem 3.3 the conclusion is immediate.\\

Now, we consider some special cases, wherein our results deduce
several well-known fixed point theorems of the existing literature.
\begin{enumerate}
\item [{(1)}]
Under the universal relation ($i.e.$ $\mathcal{R}=X^2$), Theorem 3.3
deduces to the Jotic fixed point theorem  proved in \cite{B4}, which is a
generalization of Boyd-Wong fixed point theorem \cite{B2}.
\item [{(2)}] On setting $\mathcal{R}=\preceq,$ the partial order in
Theorem 3.1 as well as Corollary 3.4, we obtain Theorem 1.3. Clearly,
$T$-closedness of $\preceq$ is equivalent to increasing property
of $T$.
\item [{(3)}] Taking $\mathcal{R}=${\it M}, the transitive binary relation in
Theorem 3.1, we obtain Theorem 1.6.
\end{enumerate}

\noindent \textbf{Conclusion:} In order to ensure the existence of fixed points for linear
contraction mapping $T$, the underlying binary relation is required to be
$T$-closed (see Theorem 2.14). But whenever, we extend Theorem 2.14 from linear contractions to
Boyd-Wong type nonlinear contractions then this restriction on the
underlying binary relation is not enough. We additionally do require
locally $T$-transitivity of $\mathcal{R}$, which
substantiates the utility of this extension. As possible problems, authors encourage the researchers of this domain to prove such results for other types of
contractions.\\

\noindent \textbf{Acknowledgements:} Both the authors are grateful to an anonymous learned referee for his fruitful comments and valuable suggestions specially pointing out an error in the earlier proof. The authors are also grateful to Professor Mihai Turinici for providing the file of his book chapter. The first author is thankful to
University Grant Commission, New Delhi, Government of India for
awarding BSR (Basic Scientific Research) Fellowship.\\

\end{document}

%% file: prembu.tex
\newcommand{\vs}{\vspace}
\newcommand{\bc}{\begin{center}}
\newcommand{\ec}{\end{center}}